\documentclass{article}
\usepackage{amssymb}
\usepackage{amscd}
\usepackage{amstext}
\usepackage{amsthm}
\usepackage{amsmath}
\usepackage{latexsym}
\title{}
\begin{document}

\theoremstyle{plain}
\newtheorem{theorem}{Theorem}[section]
\newtheorem{corollary}[theorem]{Corollary}
\newtheorem{proposition}[theorem]{Proposition}
\newtheorem{lemma}[theorem]{Lemma}
\newtheorem{definition}[theorem]{Definition}

\newtheorem{remark}{Remark}[section]
\newtheorem{example}{Example}[section]
\newtheorem{exercise}{Exercise}[section]

\title{\bf On commuting probabilities in  finite groups and rings}

\author{\bf Martin Jur\' a\v s\\
  SCAD, Savannah, GA 31401, USA\\
\bigskip
 mjuras@scad.edu, martinjuras@gmail.com\\
 \bf Mihail Ursul\\
 Department of Mathematics\\
 PNG University of Technology, Lae, PNG\\
 mihail.ursul@gmail.com
 }

\date{} \maketitle{}

\bigskip

\centerline {\bf Abstract.}
We show that the set of all commuting probabilities in finite rings is a  subset of
the set of all commuting probabilities in finite nilpotent groups of class $\leq 2$. We believe that these two sets are equal; we prove they are equal, when restricted to groups and rings with odd number of
elements.

\bigskip

\noindent {\bf Mathematics Subject Classification (2020).} Primary: 16U80. Secondary: 05C25, 20P05, 16N40, 20D15.

\bigskip

\noindent {\bf Key Words and Phrases.}  Finite group, finite ring,
commuting probability, annihilating probability, nilpotent group, nilpotent ring.

%
%

\section{Introduction and preliminaries}

In 1940, Philip Hall \cite{hall40} introduced the notion of the commuting probability in groups.
Feit and Fine \cite{feit60}, derived a combinatorial formula and a generating
function for commuting probability in matrix rings over finite fields.
In the second half of 1960's, the series of papers \cite{erdt65}, \cite{erdt67a}, \cite{erdt67b}, \cite{erdt68} by Erd\"{o}s and Tur\' an, gave birth to the statistical group theory.
In the fourth paper, among other results, the authors derived a lower bound for commuting
probability in a finite group of order $n$, and showed that the commuting probability in the symmetric group $S_n$ is asymptotically equal to $\frac{1}{n}\cdot$

A number of research and expository papers on commuting probability in groups  appeared
during late sixties and the seventies:
Joseph \cite{joseph69}, \cite{joseph77}, Galagher \cite{gala70},
Gustafson \cite{gusta73}, Machale \cite{mac74}, and Rusin \cite{rusin79}, to name a few.\footnote{
Dixon \cite{dixon02}, provides an extensive list of publications on statistical group theory in the references, up to the year 2002.}
Rusin \cite{rusin79},  characterized all finite groups  with
commuting probability $> \frac{11}{32}$.
In the nineties, Lescot \cite{lescot95}, re-derived classification of groups with commuting
probability $> \frac{1}{2}$, using the notion of isoclinism in groups introduced by Hall \cite{hall40}.

There has also been interest in the study of commuting probability of other algebraic structures; 
MacHale \cite{mac77}, investigated the notion of commuting probability in rings.
Commuting probability in semigroups has been studied in \cite{givens08},
\cite{mac90}, \cite{ponose12} and \cite{soule14}.

The dawn of the twenty-first century has seen a renewed interest
in the study of the commuting probability in groups and rings, and other types of probabilities
 in rings, such as anticommuting and annihilating probability.
Papers \cite{she00}, \cite{dixon02}, \cite{guro06}, \cite{dasna11} and \cite{hegar13}
deal with commuting probability in finite groups.
Buckley et. al.,\cite{bucmac17} and \cite{bucmaczel14}, classified all finite rings with commuting 
probability $\geq \frac{11}{32}$ and anticommuting  probability $\geq \frac{15}{32}$, respectively.

Throughout this paper,
$|A|$ denotes cardinality of the set $A.$
$Z(G)$  denotes its {\it center} of a group $G$. For $a,b\in G$,
$[a,b]=a^{-1}b^{-1}ab$ denotes the commutator of $a$ and $b$,
and $[G,G]$ denotes the {\it derived} subgroup of $G$ generated by all commutators in $G$.
Recall that $G$ is {\it nilpotent of class} $n$,
if its lower central series (of normal subgroups) terminate in the trivial
subgroup  after $n$ steps, i.e.
$$
G=G_{0}\triangleright G_{1}\triangleright \dots \triangleright G_{n}=\{e_G\},
$$
where
$G_{i}=[G_{i-1},G]$ for $i=1,2,\ldots,n$, and $G_{n-1}\neq \{e_G\}.$
{\it Commuting probability}\footnote{Some publications use the term
{\it commuting degree} in place of the {\it commuting probability}.}
in a group $G$ is defined to be the number
$$
 \text{Pr$_{c}$}(G)=\frac{|\{(a,b)\in G\times G\, : \, ab=ba\}|}{|G|^2}.
$$
For class $\mathcal G$ of finite groups, the set
$\mathfrak S_c(\mathcal G)=\{\text{Pr}_c(G) :  G\in\mathcal G\}$ is called
the {\it commuting spectrum} of $\mathcal G$.

\medskip

Rings are not assumed to be associative or unitary.  By $R(+)$ we denote the additive group of $R$.

Recall that a ring  $R$ is called {\it antisymmetric} if for all $a,b \in R$, $ab=-ba.$
$R$ is called {\it strongly antisymmetric} if the {\it dinipotent condition}, $a^2=0$,
is satisfied for all $a \in R.$
Strong antisymmetry implies antisymmetry.
A ring $R$ is said to be {\it nilpotent class $\leq n$} if the product
of any $n$ elements with any correct distribution of brackets is zero.
For a prime $p$, $R$ is called a   $p$-ring if
$|R|=p^n$ for some positive integer $n.$

The symbol $[\cdot,\cdot]$ denotes the commutator in both a group $G$ and a   ring $R$;
whenever needed, we will write $[\cdot,\cdot]_G$ and $[\cdot,\cdot]_R$
to distinguish between the two cases.

Buckley \cite{buc15}, introduced the following generalization of the notion of commuting
probability in rings.
Let $f(X,Y)=aXY+bYX$ be a formal "non-commutative polynomial" with integer
coefficients. For any   ring $R$  define a function $f^R:R\times R\rightarrow R$,
$(x,y)=axy+byx$. Let
$$
\text{Pr}_f(R)=\frac{|\{(x,y)\in R\times R \, :\,  f^R(x,y)=0\}|}{|R|^2} \cdot
$$
For class $\mathcal R$ of finite  rings,  the set
$\mathfrak S_f(\mathcal R)=\{\text{Pr}_f(R)  :  R\in\mathcal R\}$
is called the $f$-{\it spectrum} of $\mathcal R$.

Here, we are going to be mostly concerned with
the {\it commuting spectrum}, $\mathfrak S_c(\mathcal R)$ and the
{\it annihilating spectrum}, $\mathfrak S_{ann}(\mathcal R)$, with the
associated formal  "non-commutative polynomials"  $f(X,Y)=XY-YX$ and
$f(X,Y)=XY$, respectively.
The {\it commuting probability} and the {\it annihilating probability}
in a   ring $R$
are denoted by $\text{Pr}_c(R)$ and $\text{Pr}_{ann}(R)$, respectively.

We will use  the following classes of groups and rings:

\medskip

$\mathfrak G$ the class of finite groups;

$\mathfrak{G}_{nil}$ the class of finite nilpotent groups;

$\mathfrak{G}_{nil}^{(2)}$ the class of finite nilpotent  groups of class $\leq 2$;

$\mathfrak R$ the class of finite rings;

$\mathfrak{R}_{nil}^{(2)}$ the class of finite nilpotent rings of class $\leq 3$;


$\mathfrak{R}_{sa}$ the class of finite strongly antisymmetric   rings;

$\mathfrak{R}_p $ the class of $p$-rings;

for class $\mathcal C$ of finite sets, denote
$
\textup{ODD}(\mathcal C)=\{ A\in\mathcal C \, : \, \text{$|A|$ is odd}\}.
$

\medskip

Recall the following well know construction.
For a given    ring $R$, we construct the ring $N(R)$ in the following way:
the additive group of $N(R)$ is $(R\times R)(+)$
with multiplication $(a,x)(b,y)=(0,ab).$ The following Lemma is immediate.

\begin{lemma}\label{L.1.1}%
Let $R$ be a  ring. Then $N(R)$ is a nilpotent ring  of class at most $3$.
Furthermore, if $f(X,Y)=aXY+bYX$ is a formal non-commutative polynomial with integer coefficients
and $R$ is finite, then
$$
\text{Pr}_f(R)=\text{Pr}_f(N(R)).
$$
\end{lemma}

In particular, the Lemma implies
\begin{equation}
\mathfrak S_f(\mathfrak R)= \mathfrak S_f(\mathfrak{R}_{nil}^{(2)}). \label{eq1}
\end{equation}

\medskip

Ever since it was discovered that there are no finite groups with commuting probability
in the open interval  $(1, \frac{5}{8})$,
there has been an interest to understand the structure of the commuting spectrum of groups,
and later, the structure of the commuting spectrum  of rings and semigroups. The
commuting spectrum for semigroups turned out to be the simplest to understand.
Givens \cite{givens08} showed that the commuting spectrum for semigroups is dense in the
interval $[0,1]$. Later Ponomarenko and Selinski \cite{ponose12} proved that for any rational number
in $r\in (0,1]$, there is a finite semigroup $S$ such that the commuting probability in $S$ is equal
to $r$.  Soule \cite{soule14} found a single family of semigroups that has this property.

Contrastingly, for groups, Hegarty \cite{hegar13} showed that for any limit point $l\in (\frac{2}{9}, 1]$
of $\mathfrak S_c(\mathfrak G)$, there is no increasing sequence of numbers
$\{a_n\}\subset \mathfrak S_c(\mathfrak G)$, such that $l=\lim_{n\rightarrow \infty} a_n$.

Recently, Buckley and MacHale investigated relations between the commuting spectra
of finite groups and rings. Comparing the structure of these two spectra for large probabilities,
 the authors formulated two conjectures, \cite{bucmac20}, page 9:

\medskip

Conjecture 1. $\mathfrak{S}_c(\mathfrak R) \subset\mathfrak{S}_c(\mathfrak G)$.

Conjecture 2. $\mathfrak{S}_c(\mathfrak R)=\mathfrak{S}_c(\mathfrak{G}_{nil})$
            or $\mathfrak{S}_c(\mathfrak R)=\mathfrak{S}_c(\mathfrak{G}_{nil}^{(2)})$.

\medskip

In this paper, we positively resolve the first conjecture and partially resolve the second.\footnote{The authors would like to thank Victor Bovdi for his interest in this paper.}

%
%

\section{Main results}

\begin{theorem}\label{T.2.1}
$\mathfrak{S}_c(\mathfrak R)
\subseteq\mathfrak{S}_c(\mathfrak{G}_{nil}^{(2)})\subseteq
\mathfrak{S}_{ann}(\mathfrak{R}_{sa}\cap\mathfrak R_{nil}^{(2)}).$
\end{theorem}

\medskip

In \cite{bucmacshe20}, the authors determined all values in $\mathfrak{S}_c(\mathfrak R)$
that are $\geq \frac{11}{32}.$ These are
$$
1, \frac{7}{16}, \frac{11}{27}, \frac{25}{64}, \frac{11}{32}, \quad \text{and} \quad
\frac{2^{2k}+1}{2^{2k+1}} \quad \text{for} \ k=1,2,3,\ldots.
$$
Thus, $\frac{1}{2}\not \in\mathfrak{S}_c(\mathfrak R).$ But, $\frac{1}{2}\in \mathfrak{S}_c(\mathfrak G)$,
(\cite{rusin79}, page 246), and so $\mathfrak{S}_c(\mathfrak R) \neq \mathfrak{S}_c(\mathfrak G)$. In particular, $\text{Pr}_c(S_3)=\frac{1}{2}$ (see \cite{joseph77}); $S_3$ denotes the symmetric group of order 3.
This, together with the first
inclusion of Theorem \ref{T.2.1}, positively resolves Conjecture 1. As for Conjecture 2,
the Theorem states
$\mathfrak{S}_c(\mathfrak R)\subseteq\mathfrak{S}_c(\mathfrak{G}_{nil}^{(2)})$.
Now that we know  $\mathfrak{S}_c(\mathfrak R)$ is a subset of the potentially smaller
one of the two sets, $\mathfrak{S}_c(\mathfrak{G}_{nil}^{(2)})$ and
$\mathfrak{S}_c(\mathfrak{G}_{nil})$ (it is unknown whether or not
$\mathfrak{S}_c(\mathfrak{G}_{nil}^{(2)}) =\mathfrak{S}_c(\mathfrak{G}_{nil})$), we ask the following
question: Does
\begin{equation}
\mathfrak{S}_c(\mathfrak R)=\mathfrak{S}_c(\mathfrak{G}_{nil}^{(2)}) \label{eq2}
\end{equation}
hold true? We don't know. But, Equation (\ref{eq2}) does hold true,
when restricted to finite groups and finite rings with odd number of elements.
In fact, we prove the following:

\begin{theorem}\label{T.2.2}
$$
\mathfrak{S}_c(\textup{ODD}(\mathfrak R))=
\mathfrak{S}_c(\textup{ODD}(\mathfrak{G}_{nil}^{(2)}))=
\mathfrak{S}_{ann}(\textup{ODD}(\mathfrak{R}_{sa}\cap\mathfrak{R}_{nil}^{(2)})).
$$
\end{theorem}

\medskip

Next, we would like to formulate a condition, purely in terms of probabilities in rings,
that would imply Equation (\ref{eq2}).
Using Theorem \ref{T.2.1}, one obvious choice could be
$\mathfrak{S}_{ann}(\mathfrak R_{sa}\cap\mathfrak R_{nil}^{(2)})\subseteq \mathfrak{S}_c(\mathfrak R).$
We can do slightly better. Because things are working smoothly when
restricted to rings with odd number of elements, it is sufficient to focus on the "trouble makers"
which are the 2-rings.

\begin{proposition} \label{P.2.3}
If
$\mathfrak S_{ann}(\mathfrak{R}_{sa} \cap\mathfrak R_{nil}^{(2)}
 \cap \mathfrak{R}_2)\subseteq \mathfrak S_c(\mathfrak{R})$,
then Equation $ (\ref{eq2})$ holds true.
\end{proposition}

The condition of Proposition \ref{P.2.3} implies a stronger statement: If
$\mathfrak S_{ann}(\mathfrak{R}_{sa} \cap\mathfrak R_{nil}^{(2)}
 \cap \mathfrak{R}_2)\subseteq \mathfrak S_c(\mathfrak{R})$,
then both inclusions in Theorem \ref{T.2.1} can be replaced by equal signs.
Note that if there is a counterexample to the condition above, i.e. if there exists a ring $R$ such that
$R\in \mathfrak{R}_{sa} \cap\mathfrak R_{nil}^{(2)}\cap \mathfrak{R}_2$  and
$\text{Pr}_{ann}(R)\not\in \mathfrak S_c(\mathfrak{R})$,
then $\text{Pr}_{ann}(R)< \frac{11}{32}$. We conjecture that
$\mathfrak{S}_c(\mathfrak R)=\mathfrak{S}_{ann}(\mathfrak{R}_{sa}).$

%
%

\section{Proofs}

Let $N$ be an associative  nilpotent ring of class $n$.  Then $N,$ endowed with {\it "circular
multiplication",} $a\circ b=a+b+ab,$ is a group  which we will denote by
$G_N$.\footnote{Another way to associate a group to a ring such that their
commuting probabilities equate can be obtained by modifying a construction of  Mal'cev \cite{mal65}.
For an arbitrary ring $R$, define a binary operation on
$R\times R$ by $(a,b)\cdot (c,d)=(a+c,ac+b+d)$. This operation is associative, has unit $(0,0)$ and
$(a,b)^{-1}=(-a,a^2-b).$  $G=(R\times R, \cdot\ )$ is a nilpotent group of class at most 2 and $\text{Pr}_c(R)=\text{Pr}_c(G)$.
Note that, unlike the construction of $G_N$, the ring $R$ is not required to be nilpotent or
associative!}
$0$ is the unit element in $G_N$ and
$a^{-1}=-a+a^2-a^3+\cdots +(-1)^{n-1}a^{n-1}$ is the inverse of
$a$ in $G_N,$  $a\circ a^{-1}= a^{-1}\circ a=0.$
 Since, $ab=ba$ if and only if $a\circ b=b\circ a,$ then, if  $N$ is finite,
\begin{equation}
\text{Pr}_c(N)=\text{Pr}_c(G_{N}), \label{eq3}
\end{equation}

\begin{lemma}\label{L.3.1} %
 Let $N$ be a  nilpotent ring of class at most $3$ (hence, also an associative ring).
 Let $a,b,c\in N.$ Then

 \noindent (i)  \ \ $[a,b]_{G_N}=[a,b]_N,$

 \noindent (ii) \ $[a,b]_{G_N}\circ [c,d]_{G_N}=[a,b]_N +[c,d]_N,$

 \noindent (iii) $G_N$ is a nilpotent group of class $\leq 2$.
\end{lemma}
\proof%
$(i)$  follows by direct computation.

\noindent $(ii)$. By $(i)$,
$$
[a,b]_{G_N}\circ [c,d]_{G_N}=[a,b]_N \circ [c,d]_N
$$
$$
= [a,b]_N +[c,d]_N+[a,b]_N[c,d]_N =[a,b]_N +[c,d]_N.
$$
$(iii)$. By (i), $[[a,b]_{G_N},c]_{G_N}= [[a,b]_{N},c]_{N}=0$.
\hfill $\rule{0.5em}{0.5em}$

\medskip

Let $G$ be a nilpotent group of class $\leq 2$ and let $Z=Z(G)$ be the center of $G$.
Then $G/Z$ is abelian.  By $R_G$,
denote the ring with the additive group $G/Z\oplus Z$, and the
multiplication defined by
\begin{equation}
(aZ,x)\cdot (bZ,y)=(Z,[a,b]), \label{eq4}
\end{equation}
where $[a,b]=a^{-1}b^{-1}ab$ \, is the commutator in $G$.
Explicitly, the addition in $R_G$ is given by
$$
(aZ,x)+(bZ,y)=(abZ,xy).
$$
$(Z,  e_G)$ is the zero element and $(a^{-1}Z,x^{-1})$ is the additive inverse
of $(aZ,x)$.

To verify that $R_G$ is indeed a   ring,
the distributive laws have to be satisfied.
Let $a,b,c\in G$ and $x,y,z\in Z$.
We have
$$
(cZ,z)\cdot ((aZ,x)+(bZ,y))=(cZ,z)\cdot (abZ,xy)
=(Z,[c,ab]).
$$
On the other hand,
$$
(cZ,z)\cdot (aZ,x)+(cZ,z)\cdot (bZ,y))
=(Z,[c,a])+(Z,[c,b])
=(Z,[c,a][c,b]).
$$
Using $[G,G]\subseteq Z$, we deduce
$$
[c,a][c,b]=c^{-1}a^{-1}cac^{-1}b^{-1}cb=c^{-1}[a,c^{-1}]b^{-1}cb=c^{-1}b^{-1}[a,c^{-1}]cb
$$
$$
=c^{-1}b^{-1}a^{-1}cac^{-1}cb=c^{-1}b^{-1}a^{-1}cab=c^{-1}(ab)^{-1}c(ab)=[c,ab].
$$

Hence, the left distributive law is satisfied. The proof of the right distributive law is
 similar.\footnote{Proposition 3 \cite{bucmac20}, states that the condition $[c,a][c,b]=[c,ab]$
 for all $a,b, c\in G$ is equivalent to $G$ being nilpotent of class $\leq 2$.}

\medskip

\begin{lemma}\label{L.3.2}
Let $G$ be a nilpotent group of class at most $2$. Then $R_G$ is a strongly antisymmetric
nilpotent ring of class at most $3$.
If $G$ is finite, then $|R_G|=|G|$ and
\begin{equation}
\textup{Pr}_c(G)=\textup{Pr}_{ann}(R_G). \label{eq5}
\end{equation}
\end{lemma}
\proof%
$|R_G|=|G/Z||Z|=|G|$.
$R_G^3=0$ and strong antisymmetry of $R_G$ follows immediately from the multiplication
formula (\ref{eq4}) and the fact that $G$ is a nilpotent group of class $\leq 2$.
To prove (\ref{eq5}), it suffices to note that\\
$
(aZ,x)\cdot (bZ,y)=(Z,[a,b])=(Z, e_G)
$
if and only if
$[a,b]=e_G.
$ But, this is exactly when $ab=ba.$
\hfill $\rule{0.5em}{0.5em}$

\medskip

\noindent {\bf Proof of Theorem \ref{T.2.1}.}

We first show that
$\mathfrak{S}_c(\mathfrak R)\subseteq\mathfrak{S}_c(\mathfrak{G}_{nil}^{(2)}).$
Let $r\in \mathfrak{S}_c(\mathfrak R)$ By Lemma \ref{L.1.1}, there is a nilpotent
ring $N$ of class at most 3 such that $r=\text{Pr}_c(N)$.
By Lemma \ref{L.3.1}$(iii)$, $G_{N}$ is a nilpotent group of class at most 2 and
by Equation (\ref{eq3}), $\text{Pr}_c(G_N)=\text{Pr}_c(N)$. We conclude that
$r\in\mathfrak{S}_c(\mathfrak{G}_{nil}^{(2)}).$

To prove the second inclusion,
consider $G\in \mathfrak{G}_{nil}^{(2)}$. By Lemma \ref{L.3.2}

\noindent $R_G\in \mathfrak R_{sa}\cap\mathfrak{R}_{nil}^{(2)}$ \,  and
$\text{Pr}_{ann}(R_G)=\text{Pr}_c(G).$
\hfill $\rule{0.5em}{0.5em}$

\medskip

\begin{lemma}\label{L.3.3}
Let $R$ be a finite  antisymmetric ring and with odd number of elements. Then
$$
\textup{Pr}_c(R)=\textup{Pr}_{ann}(R).
$$
\end{lemma}
\proof%
In an antisymmetric   ring, $ab=-ba$. Hence, $ab=ba$ iff $2ab=0$.
Since $|R|$ is odd, $2ab=0$ iff $ab=0$.
\hfill $\rule{0.5em}{0.5em}$

\medskip

\noindent {\bf Proof of Theorem \ref{T.2.2}.}

To prove
$
\mathfrak{S}_c(\textup{ODD}(\mathfrak R))\subseteq
\mathfrak{S}_c(\textup{ODD}(\mathfrak{G}_{nil}^{(2)}))\subseteq
\mathfrak{S}_{ann}(\textup{ODD}(\mathfrak{R}_{sa}\cap\mathfrak{R}_{nil}^{(2)})),
$
we follow the proof of Theorem  \ref{T.2.1}
and note that $|N|=|G_N|$ and $|G|=|R_G|$.

To conclude the proof of Theorem \ref{T.2.2}, it suffices to show
$$
\mathfrak{S}_{ann}(\textup{ODD}(\mathfrak{R}_{sa}\cap\mathfrak{R}_{nil}^{(2)}))\subseteq
\mathfrak{S}_c(\textup{ODD}(\mathfrak{R})).
$$
Let
$r \in \mathfrak{S}_{ann}(\textup{ODD}(\mathfrak{R}_{sa}
\cap\mathfrak{R}_{nil}^{(2)}))$ and let $R \in \textup{ODD}(\mathfrak{R}_{sa}\cap\mathfrak{R}_{nil}^{(2)})$
such that $r=\text{Pr}_{ann}(R).$ By Lemma \ref{L.3.3},
$r=\text{Pr}_{ann}(R)=\text{Pr}_{c}(R) \in \mathfrak{S}_c(\textup{ODD}(\mathfrak{R})).$
\hfill $\rule{0.5em}{0.5em}$

\medskip

In a   ring, the additive order of $ab$ divides the additive orders of both
$a$ and $b$. In particular, if the additive orders of $a$ and $b$
are relatively prime, then $ab=0$. As a consequence of this fact and the
GH fundamental theorem  of finite abelian groups, GH a finite   ring
is a product of   $p$-rings. In turn, this implies that
$$
\text{Pr}_{ann}(R_1\times R_2)= \text{Pr}_{ann}(R_1)\, \text{Pr}_{ann}(R_2),
$$
for any two rings $R_1$, $R_2$.\footnote{Similar arguments show that
$\text{Pr}_{f}(R_1\times R_2)= \text{Pr}_{f}(R_1)\, \text{Pr}_{f}(R_2)$, for
any noncommutative formal polynomial $f(X,Y)=aXY+ bYX.$}

 We say that a class  $\mathcal R$ of finite rings is {\it hereditary},
if any  subring of a ring in $\mathcal R$ is also in $\mathcal R$.

Let $p$ be a prime number and
assume that  a class $\mathcal C$ of finite rings is hereditary.
Then $\mathcal C_p= \mathcal C\cap \mathfrak{R}_p$ also is hereditary.
Furthermore,
$\mathfrak S_{ann}(\mathcal C)$ and $\mathfrak S_{ann}(\mathcal C_p)$,
are monoids and values in
$\mathfrak S_{ann}(\mathcal C)$ are finite products of values taken from the set
$\bigcup_p\mathfrak S_{ann}(\mathcal{C}_p)$,
where $p$ runs  all prime number.

\medskip

\noindent {\bf Proof of Proposition \ref{P.2.3}.}

It is easy to see that the class
$\mathcal C= \mathfrak{R}_{sa} \cap\mathfrak R_{nil}^{(2)}$
is hereditary.
Assume
$\mathfrak S_{ann}(\mathcal C \cap \mathfrak{R}_2)\subseteq \mathfrak S_c(\mathfrak{R})$.
If $p\neq 2$ is a prime, by Theorem \ref{T.2.2},
$\mathfrak{S}_{ann}(\textup{ODD}(\mathcal C))
=\mathfrak{S}_c(\textup{ODD}(\mathfrak R))$ and so
$\mathfrak S_{ann}(\mathcal C\cap\mathfrak{R}_p)
\subseteq \mathfrak S_c(\mathfrak R)$.
Hence, for all primes  $p$, the monoids
$\mathfrak S_{ann}(\mathcal C\cap\mathfrak{R}_p)
\subseteq \mathfrak S_c(\mathfrak R)$.
Using the considerations above, we conclude
$\mathfrak S_{ann}(\mathcal C)\subseteq\mathfrak S_c(\mathfrak{R})$.
By Theorem \ref{T.2.1}, the reverse inclusion is satisfied, and so
the proposition follows.
\hfill $\rule{0.5em}{0.5em}$


\end{document}